\documentclass[a4paper,12pt,final]{amsart}
\usepackage{times,a4wide,mathrsfs,amssymb,amsmath,amsthm,enumerate,xypic,tikzsymbols,dsfont}

\newcommand{\C}{\mathbb{C}}
\newcommand{\ZZ}{\mathbb{Z}}

\newcommand{\QQ}{\mathbb{Q}}
\newcommand{\NN}{\mathbb{N}}
\newcommand{\PP}{\mathbb{P}}

\newcommand{\OO}{\mathcal O}

\newcommand{\Sy}{\mathfrak S}

\newcommand{\A}{\operatorname{A}}
\newcommand{\DA}{\operatorname{DA}}
\newcommand{\N}{\operatorname{N}}

\newcommand{\MM}{\mathcal M}

\newcommand{\wt}{\widetilde}
\newcommand{\rom}{\romannumeral}

\newcommand{\one}{\mathds{1}}

\DeclareMathOperator{\aut}{Aut}
\DeclareMathOperator{\bir}{Bir}

\DeclareMathOperator{\sym}{Sym}

\DeclareMathOperator{\jac}{Jac}

\newtheorem{theorem}{Theorem}[section]

\newtheorem{lemma}[theorem]{Lemma}

\newtheorem{proposition}[theorem]{Proposition}
\newtheorem{conjecture}[theorem]{Conjecture}
\newtheorem{convention}{Conventions}

\newtheorem{nonumbering}{Theorem}

\theoremstyle{definition}
\newtheorem{remark}[theorem]{Remark}
\newtheorem{definition}[theorem]{Definition}
\newtheorem{question}[theorem]{Question}

\newtheorem{nonumberingt}{Acknowledgments}

\begin{document}

\author[Robert Laterveer]
{Robert Laterveer}

\address{Institut de Recherche Math\'ematique Avanc\'ee,
CNRS -- Universit\'e 
de Strasbourg,\
7 Rue Ren\'e Des\-car\-tes, 67084 Strasbourg CEDEX,
FRANCE.}
\email{robert.laterveer@math.unistra.fr}

\title{On the Chow ring of some special Calabi--Yau varieties}

\begin{abstract} We consider Calabi--Yau $n$-folds $X$ arising from certain hyperplane arrangements. Using Fu--Vial's theory of distinguished cycles for varieties with motive of abelian type, we show that the subring of the Chow ring of $X$ generated by divisors, Chern classes and intersections of subvarieties of positive codimension injects into cohomology. We also prove Voisin's conjecture for $X$, and Voevodsky's smash-nilpotence conjecture for odd-dimensional $X$.
\end{abstract}


\thanks{\textit{2020 Mathematics Subject Classification:}  14C15, 14C25, 14C30}
\keywords{Algebraic cycles, Chow group, motive, Bloch--Beilinson filtration, distinguished cycles, section property}
\thanks{Supported by ANR grant ANR-20-CE40-0023.}

\maketitle

\section{Introduction}

Given a smooth projective variety $Y$ over $\C$, let $\A^i(Y):=\operatorname{CH}^i(Y)_{\QQ}$ denote the Chow groups of $Y$ (i.e. the groups of codimension $i$ algebraic cycles on $Y$ with $\QQ$-coefficients, modulo rational equivalence \cite{F}, \cite{MNP}, \cite{Vo}). The intersection product defines a ring structure on $\A^\ast(Y)=\bigoplus_i \A^i(Y)$, the {\em Chow ring\/} of $Y$.
In the case of K3 surfaces, this ring structure has a remarkable property:

\begin{theorem}[Beauville--Voisin \cite{BV}]\label{K3} Let $S$ be a projective K3 surface. 
The $\QQ$-subalgebra
  \[  R^\ast(S):=  \bigl\langle  \A^1(S), c_j(S) \bigr\rangle\ \ \ \subset\ \A^\ast(S) \]
  injects into cohomology under the cycle class map.
  \end{theorem}

A fancy way of rephrasing this result is as follows: for any variety $Y$, let 
  \[ \N^i(Y):= \A^i(Y)/\A^i_{hom}(Y)\] 
  denote the quotient, where $\A^i_{hom}(Y)\subset \A^i(Y)$ denotes the homologically trivial cycles. Then for K3 surfaces $S$ (which have $H^1(S,\QQ)=0$ and so $A^1(S)$ injects into cohomology), Theorem \ref{K3} says that the
$\QQ$-algebra epimorphism 
  \[ \A^\ast(S)\ \twoheadrightarrow\ \N^\ast(S)\] 
  admits a section, whose image contains the Chern classes of $S$ -- that is, $S$ has the {\em section property\/}, in the language of \cite{FV}. 

It is then natural to ask which other varieties have the section property. An interesting partial answer is given in \cite{FV}, by extending O'Sullivan's theory of distinguished cycles from abelian varieties to varieties with motive of abelian type: if a variety $Y$ verifies condition $(\star)$ of loc. cit., then all powers of $Y$ have the section property (cf. subsection \ref{ss:fv} below). One could say that varieties verifying the condition $(\star)$ form a kind of ``meilleur des mondes possibles'', a world in which Chow motives and their multiplicative behaviour are well-understood. Unfortunately, inhabitants of this meilleur des mondes are rather scarce;
some examples of  varieties verifying condition $(\star)$ are given in \cite{FV} and \cite{LV}. 

The main result of the present paper exhibits special Calabi--Yau varieties of any dimension for which the Chow ring is just as well-behaved as that of K3 surfaces:

\begin{nonumbering}[=Theorem \ref{main}] Let $X$ be a hyperelliptic Calabi--Yau variety of dimension $n\ge 2$.
The $\QQ$-subalgebra 
  \[  R^\ast(X):=   \Bigl\langle  \A^1({X}), \, \A^i({X})\cdot \A^j({X}),\, c_k(X)\Bigr\rangle\ \ \ \subset\ \A^\ast(X)\ \ \ \ \ (i,j>0) \]
generated by divisors, Chern classes and by cycles that are intersections of two cycles of positive codimension injects into cohomology. In particular, the image of the intersection product
  \[ \A^i(X)\otimes \A^{j}(X)\ \ \to\ \A^{i+j}(X)\ \ \ \ (i,j>0) \]
  injects into cohomology.
  \end{nonumbering}
  
 Here, a {\em hyperelliptic\/} Calabi--Yau variety is defined as follows: given $2n+2$ hyperplanes in general position in $\PP^n$, the double cover of $\PP^n$ branched along the union of hyperplanes admits a crepant resolution that is Calabi--Yau. We say that the resulting Calabi--Yau $n$-fold is {\em hyperelliptic\/} if the hyperplanes osculate a rational normal curve (cf. subsection \ref{ss:he} below). These Calabi--Yau varieties have been studied in \cite{MT}, \cite{GSSZ}, \cite{Tera2}.
 
 The behaviour exhibited by Theorem \ref{main} is remarkable, in the sense that for general Calabi--Yau varieties $Y$ one does {\em not\/} expect that the image of the intersection product
  \[ \A^i(Y)\otimes \A^{j}(Y)\ \ \to\ \A^{i+j}(Y)\ \ \ \ (i,j>0) \]
  injects into cohomology; one only expects this for $i+j=\dim Y$ (and this last expectation is known for complete intersection Calabi--Yau varieties \cite{Fu}, but wide open in general).
  
  In proving Theorem \ref{main}, we rely on the ``meilleur des mondes'' formalism of \cite{FV}. We actually prove that a certain blow-up of $X$ verifies condition $(\star)$ of loc. cit.; the result for $X$ is then a consequence of the nice behaviour of the formalism. This reasoning is very similar to that of \cite{LV}.
 
 As a by-product of the argument,
 we also obtain some new cases where Voisin's conjecture \cite{V9} is verified:
 
 \begin{nonumbering}[=Theorem \ref{conjV}] Let $X$ be a hyperelliptic Calabi--Yau $n$-fold. Any two zero-cycles $a,a^\prime\in \A^n_{hom}(X)$ satisfy
		\[ a\times a^\prime = (-1)^n \, a^\prime\times a\ \ \ \hbox{in}\ \A^{2n}(X\times
		X)\ .\]
		(Here, $a\times a^\prime$ denotes the exterior product $(p_1)^\ast(a)\cdot
		(p_2)^\ast(a^\prime)\in \A^{2n}(X\times X)$, where $p_j$ is projection to the
		$j$-th factor.)
	\end{nonumbering}	
	
Another by-product concerns a conjecture of Voevodsky \cite{Voe}:

\begin{nonumbering}[=Theorem \ref{conjVoe}] Let $X$ be a hyperelliptic Calabi--Yau variety of odd dimension. Then homological equivalence and smash-equivalence coincide for all algebraic cycles on $X$.
\end{nonumbering}

 The aim of this paper is twofold: on the one hand, we want to promote the ``meilleur des mondes'' formalism of \cite{FV} (and encourage others to find new instances where this formalism can be applied); on the other hand, we want to raise interest for questions concerning the multiplicative structure of the Chow ring of varieties (and to this end, we have included some open questions concerning other Calabi--Yau varieties, cf. section \ref{quest}).

 \vskip0.6cm

\begin{convention} In this article, the word {\sl variety\/} will refer to a reduced irreducible scheme of finite type over $\C$. A {\sl subvariety\/} is a (possibly reducible) reduced subscheme which is equidimensional. 

{\bf All Chow groups are with rational coefficients}: we will denote by $\A_j(Y)$ the Chow group of $j$-dimensional cycles on $Y$ with $\QQ$-coefficients; for $Y$ smooth of dimension $n$ the notations $\A_j(Y)$ and $\A^{n-j}(Y)$ are used interchangeably. 
The notations $\A^j_{hom}(Y)$ and $\A^j_{AJ}(Y)$ will be used to indicate the subgroup of homologically trivial (resp. Abel--Jacobi trivial) cycles.
For a morphism $f\colon X\to Y$, we will write $\Gamma_f\in \A_\ast(X\times Y)$ for the graph of $f$.

The covariant category of Chow motives (i.e., pure motives with respect to rational equivalence with $\QQ$-coefficients as in \cite{An0})
will be denoted $\MM_{\rm rat}$.
\end{convention}

\section{Preliminaries}

\subsection{Intersection theory on quotient varieties}

\begin{lemma}\label{quotient} Let $M$ be a {\em quotient variety\/}, \emph{i.e.}
		$M=M^\prime/G$ where $M^\prime$ is a smooth quasi-projective variety and
		$G\subset\aut(M^\prime)$ is a finite group.
		Then $\A^\ast(M):=\oplus_i \A_{\dim M-i}(M)$ is a commutative
		graded ring, with the usual functorial properties.
	\end{lemma}
	
	\begin{proof} According to \cite[Example 17.4.10]{F}, the natural map
		\[ \A^i(M)\ \to\ \A_{\dim M-i}(M) \]
		from operational Chow cohomology (with $\QQ$-coefficients) to the usual Chow
		groups (with $\QQ$-coefficients) is an isomorphism. The lemma follows from the
		good formal properties of operational Chow cohomology.
	\end{proof}

\begin{remark} In particular, Lemma \ref{quotient} implies that the formalism of correspondences and pure motives (with $\QQ$-coefficients) makes sense for projective quotient varieties.
\end{remark}

\subsection{Hyperelliptic Calabi--Yau varieties}
\label{ss:he}

It is a well-known fact that hyperplane arrangements give rise to Calabi--Yau varieties:

\begin{proposition}\label{p1} Let $H_1,\ldots,H_{2n+2}$ be hyperplanes in $\PP^n$ that are in general position (i.e. $\dim H_{i_1}\cap\cdots \cap H_{i_j} = n-j$ for each subset $\{i_1,\ldots,i_j\}\subset\{1,\ldots,2n+2\}$). Let $\bar{X}\to\PP^n$ be the double cover ramified along $\cup_{i=1}^{2n+2} H_i$. Then $\bar{X}$ is a quotient variety, and there exists a resolution of singularities $f\colon X\to\bar{X}$ such that $X$ is a Calabi--Yau variety. The morphism $f$ is a sequence of blow-ups with smooth centers $Z_i$ that have trivial Chow groups (i.e. $\A^\ast_{hom}(Z_i)=0$).
\end{proposition}

\begin{proof} As explained in \cite{GSSZ}, the arrangement $\{H_i\}$ is determined by an $(n+1)\times(2n+2)$-matrix $(b_{ij})$, where $b_{ij}\in\C$.
Let $Y\subset\PP^{2n+1}$ be the complete intersection of quadrics
   \[   \begin{cases} b_{00} x_0^2+  b_{01}x_1^2+\cdots \  \cdots \ \cdots+ b_{0,2n+1}x_{2n+1}^2&=0\ ,\\
                          \ \ \ \vdots\\
                         \ \ \  \vdots\\
                           b_{n0} x_0^2+b_{n1} x_1^2+\cdots \cdots \cdots + b_{n,2n+1} x_{2n+1}^2&=0\ .\\     
                           \end{cases}\]
       The non-singularity of $Y$ is equivalent to the $H_i$ being in general position \cite[Proposition 3.1.2]{Tera}. There is an isomorphism
       \[ \bar{X}\cong Y/G \]
       for some finite group $G$ (this is proven for $n=3$ in \cite[Proposition 2.5]{GSSZ}; the argument works for general $n$), and so $\bar{X}$ is a quotient variety.

 A crepant resolution $X\to\bar{X}$ is constructed in \cite[Section 5.1]{CH} (an alternative construction is given in \cite[Proposition 4.2]{MT}). 
For later use, we give a precise description of the resolution algorithm. The resolution $X\to\bar{X}$ is constructed as a cartesian diagram
\[  \begin{array}[c]{cccccccc}      X=:X_{m}\ \ \ \   & \xrightarrow{}\ \ \ &      \cdots & \cdots & \xrightarrow{}& X_1 & \xrightarrow{}& \ \ X_0:= \bar{X} \\       
                       &&&&&&&\\
                  \ \   \ \  \ \  \ \ \ \ \downarrow{\scriptstyle \pi_{m}} \ \ \ &&&&& \ \ \ \  \downarrow{\scriptstyle \pi_1} &&   \downarrow{\scriptstyle \pi_0}\ \ \ \ \ \ \\
                      &&&&&&&\\
                 \  \  \ \ \ \  \  P_{m}\ \ \ \   & \xrightarrow{r_{m-1}}\ \ \ &      \cdots & \cdots & \xrightarrow{}& P_1 & \xrightarrow{r_0}& \ \ P_0:= \PP^n \ ,\\
                 \end{array}\]   
where $r_0$ is a blow-up with center the codimension 2 intersection $Q_0:=H_1\cap H_2$, and
each $r_j$ is a blow-up with center $Q_j\subset P_j$, where $Q_j$ is the strict transform of an intersection $H_{i_1}\cap H_{i_2}$. 
(This description is perhaps not immediately apparent from reading \cite[Section 5.1]{CH}, but this becomes crystal clear from the reinterpretation of \cite{CH} given in \cite{IL}: the arrangement $\{H_i\}$ is {\em splayed\/}, in the sense of loc. cit. (cf. \cite[Lemmas 3.24 and 3.25]{IL}), and hence the resolution algorithm \cite[Algorithm 1.5]{IL} consists of blowing-up all pairwise intersections $H_{i_1}\cap H_{i_2}$, in arbitrary order.)
\end{proof}

\begin{remark} In case $n=2$, $X$ as in Proposition \ref{p1} is a K3 surface of the type studied in \cite{Par}, \cite{MS}, \cite{Yos}. In case $n=3$, Calabi--Yau varieties $X$ as in Proposition \ref{p1} are special cases of so-called  ``double octics''; these special cases have been intensively studied, particularly their modular properties \cite{Mey}, \cite{CK}, \cite{CM}, \cite{CM2}, \cite{CSS}, \cite{CS}, \cite{GSSZ}, \cite{Tera2}.
\end{remark}

In order to define hyperelliptic Calabi--Yau varieties, we consider a special case of the above construction:

\begin{proposition}\label{p2} Let $p_1,\ldots,p_{2n+2}\in\PP^1$ be distinct points, and let
  \[  H_i:= \gamma\bigl(  p_i\times (\PP^1)^{n-1}\bigr)\ \ \subset\ \PP^n\ ,\]
  where $\gamma\colon (\PP^1)^n\to \PP^n$ is the natural map
  \[  \gamma\colon\ \ (\PP^1)^n\to  \sym^n (\PP^1)\cong \PP H^0(\PP^1,\OO_{\PP^1}(n))\cong \PP^n  \ .\]
  Then $H_1,\ldots,H_{2n+2}$ is a hyperplane arrangement in general position, and hence gives rise to a Calabi--Yau $n$-fold $X$.
  
  The hyperplanes $H_1,\ldots,H_{2n+2}$ are tangent to a rational normal curve of degree $n$. Conversely, any hyperplane arrangement osculating a rational normal curve of degree $n$ arises in this way.
   \end{proposition}

\begin{proof} The hyperplane $H_i$ corresponds to degree $n$ divisors on $\PP^1$ containing the point $p_i$. The intersection of hyperplanes $H_{i_1}\cap\cdots \cap H_{i_r}$ then corresponds to degree $n$ divisors on $\PP^1$ containing $r$ distinct points; this has dimension equal to the expected dimension $n-r$.

The second statement (which is observed in \cite[Remark 2.10]{GSSZ}) follows from the fact that the image of the diagonal embedding of $\PP^1$ under $\gamma$ is a rational normal curve tangent to the $H_i$.
\end{proof}

\begin{definition}\label{def} Let $\{H_i\}$ be a hyperplane arrangement as in Proposition \ref{p2}, and let $\bar{X}\to\PP^n$ be the double cover branched along $\cup_i H_i$. The crepant resolution $X\to\bar{X}$ constructed via the algorithm of Proposition \ref{p1} is called a {\em hyperelliptic\/} Calabi--Yau variety.
\end{definition}

\begin{remark} The appellation ``hyperelliptic'' comes from \cite{GSSZ}, where the moduli space of these Calabi--Yau varieties of dimension $n=3$ is called the ``hyperelliptic locus''.
As observed in \cite[Remark 2.10]{GSSZ}, the Calabi--Yau $n$-folds of Proposition \ref{p1} have moduli dimension $n^2$, while the hyperelliptic Calabi--Yau $n$-folds of Proposition \ref{p2} have moduli dimension $2n-1$.

In case $n=2$, this is the well-known fact that K3 surfaces coming from double planes branched along 6 lines form a $4$-dimensional family, while imposing that the 6 lines are tangent to a conic one obtains exactly the quartic Kummer K3 surfaces (which form a 3-dimensional family). The moduli space of these surfaces is studied in \cite{MS}.
\end{remark}

The next result provides justification for Definition \ref{def}:

\begin{proposition}\label{motive} Let $\bar{X}\to\PP^n$ be a double cover branched along a hyperplane arrangement as in Proposition \ref{p2}, and let
$X\to \bar{X}$ be the crepant resolution  coming from Proposition \ref{p1} (i.e. $X$ is a hyperelliptic Calabi--Yau $n$-fold). 

\noindent
(\rom1) There is an isomorphism $\bar{X}\cong C^n/G$, where $C$ is a hyperelliptic curve of genus $n$, and $G\subset\aut(C^n)$ is a finite group of automorphisms; in particular,
 $\bar{X}$ has only quotient singularities.
 
 \noindent
 (\rom2) There exist a motive $M\in\MM_{\rm rat}$ and isomorphisms of Chow motives
  \[ \begin{split}   h(\bar{X}) &\cong   M\oplus \bigoplus \one(\ast)\ ,\\
                              h({X}) &\cong     M\oplus \bigoplus \one(\ast)\ \ \ \ \ \hbox{in}\ \MM_{\rm rat}\ ,\\
                              \end{split}\]
         where $H^j(M,\QQ)=0$ for $j\not=n$.      
 \end{proposition}

\begin{proof} 

\noindent
(\rom1) In case $n=3$, this is contained in \cite[Section 2.3]{GSSZ}; the same argument works for arbitrary $n$, as we now explain. Let $q\colon C\to\PP^1$ be the hyperelliptic curve branched along
the points $p_1,\ldots, p_{2n+2}\in\PP^1$. The morphism
   \[ h\colon\ \ C^n\ \xrightarrow{q^n}\ (\PP^1)^n\ \xrightarrow{\gamma}\ \PP^n \]
   is a Galois covering of degree $2^n\cdot n!$, with Galois group $G_2\cong \langle \iota_1,\ldots,\iota_n\rangle\rtimes \Sy_n$ (here the $\iota_j$ denotes the hyperelliptic involution of the $j$th copy of $C$). Let $G\subset G_2$ be the index 2 subgroup $G:=N\rtimes\Sy_n$, where $N$ is the kernel of the sum homomorphism $(\ZZ/2\ZZ)^n\to \ZZ/2\ZZ$. There is a cartesian diagram
     \[  \begin{array}[c]{ccc}
          C^n &\xrightarrow{p} & C^n/G\\
          &&\\
        \ \   \downarrow{\scriptstyle q^n}&&  \    \downarrow{\scriptstyle \pi}\\
          &&\\
          (\PP^1)^n &\xrightarrow{\gamma}& \PP^n\ .\\
          \end{array}\]
   The Galois group of $\pi$ is isomorphic to $G_2/G$, and is generated by the image of $\iota_1$. Hence, the ramification locus of $\pi$ is the image of the fixed locus $L_1\subset C^n$ of $\iota_1$ under $\pi\circ  p$. This is the same as the image of $L_1$ under $\gamma\circ q^n$. Since
   \[   q^n(L_1)=\bigcup_{i=1}^{2n+2} p_i\times (\PP^1)^{n-1}\ ,\] 
   it follows that the ramification locus of $\pi$ consists of the union of the hyperplanes
   \[      H_i:= \gamma\bigl(  p_i\times (\PP^1)^{n-1}\bigr)\ \ \subset\ \PP^n\ .\]
   As $C^n/G$ and $\bar{X}$ are double covers of $\PP^n$ with the same ramification locus, this proves (\rom1).
   
   As for (\rom2), note that $\bar{X}$ (being a double cover of $\PP^n$ branched along a degree $2n+2$ divisor) is isomorphic to an ample hypersurface in weighted projective space $\PP(1^{n+1},n+1)$. By weak Lefschetz, plus the fact that $\bar{X}$ is a quotient variety and hence satisfies Poincar\'e duality with $\QQ$-coefficients (cf. for instance \cite[Section 4.2.2]{Dol}), it follows that 
      \[    H^j(\bar{X},\QQ) =\begin{cases} \QQ&\hbox{\ \  if\ $j\not=n$\ even}\ ,\\
                                                   0 &\hbox{\ \  if\ $j\not=n$\ odd}\ .\\
                                                   \end{cases}\]
                                             Let $M\subset h(\bar{X})$ be the submotive defined by the projector
                                 \[ \Delta_{\bar{X}}-\sum_{j=0}^n  {1\over d}\, h^j\times h^{n-j} \ \ \ \in\ \A^n(\bar{X}\times\bar{X})      \]
                                 (where $d$ is the degree of $\bar{X}$). This gives a decomposition of $h(\bar{X})$ as requested.
                 
                 The resolution $X\to\bar{X}$ is done by blowing up subvarieties with trivial Chow groups, and so the blow-up formula gives an isomorphism
                 \[ h(X)\cong h(\bar{X})\oplus \bigoplus\one(\ast)\ \ \ \hbox{in}\ \MM_{\rm rat}\ .\]
                 This gives the requested decomposition of $h(X)$.            
%
%
%
  \end{proof}

\begin{remark} It seems likely that in Proposition \ref{motive}(\rom2), one actually has an isomorphism 
  \[ M\cong \sym^n h^1(C)\] 
  (this is stated for $n=3$ in \cite[Equation (1.3)]{Tera2} and \cite[Proposition 2.9]{GSSZ}\footnote{More precisely, it is stated in \cite{Tera2} and \cite{GSSZ} that $H^3(\bar{X},\QQ)\cong \wedge^3 H^1(C,\QQ)$ (from which the isomorphism of motives would readily follow). However the reference \cite{Tera2} contains no proof, and the proof of \cite[Proposition 2.9]{GSSZ} contains a gap: it is asserted in loc. cit. that $H^3( \sym^3(C),\QQ)\cong H^3(\jac(C),\QQ)$ where $C$ is a genus 3 curve and $\jac(C)$ its Jacobian. However, $\dim H^3( \sym^3(C),\QQ)=26$ whereas $\dim H^3(\jac(C),\QQ)=20$.}).
\end{remark}

\subsection{The section property and distinguished cycles}
\label{ss:fv}

The following notion was introduced by O'Sullivan \cite{OS}.
	
	\begin{definition}[Symmetrically distinguished cycles on abelian varieties
		\cite{OS}]\label{def:SD}
		Let $B$ be an abelian variety and $\alpha\in \A^*(B)$. For each integer
		$m\geq
		0$, denote by $V_{m}(\alpha)$ the $\QQ$-vector subspace of $\A^*(B^{m})$
		generated
		by elements of the form
		$$p_{*}(\alpha^{r_{1}}\times \alpha^{r_{2}}\times \cdots\times
		\alpha^{r_{n}}),$$
		where $n\leq m$, $r_{j}\geq 0 $ are integers, and $p :  B^{n}\to B^{m}$ is a
		closed immersion with each component $B^{n}\to B$ being either a projection or
		the composite of a projection with $[-1]: B\to B$. Then $\alpha$ is
		\emph{symmetrically distinguished} if for every $m$ the restriction of the
		projection $\A^*(B^{m})\to \N^*(B^{m})$ to $V_{m}(\alpha)$ is
		injective.
	\end{definition}
	
	The main result of \cite{OS} is\,:
	
	\begin{theorem}[O'Sullivan \cite{OS}]\label{thm:SD}
		Let $B$ be an abelian variety. Then $\operatorname{DA}^*(B)$, the
		symmetrically distinguished cycles in
		$\A^*(B)$, form a graded $\QQ$-subalgebra that contains symmetric divisors
		and that is
		stable under pull-backs and push-forwards along homomorphisms of abelian
		varieties. Moreover
		the composition $$\operatorname{DA}^*(B)\hookrightarrow
		\A^*(B)\twoheadrightarrow 
		\N^*(B)$$ is an isomorphism of $\QQ$-algebras.
	\end{theorem}
	
	Let $X$ be a smooth projective variety such that its Chow motive $h(X)$
	belongs
	to the strictly full and thick subcategory of Chow motives generated by the
	motives of abelian varieties. We say that $X$ has motive \emph{of abelian
		type}.
	A \emph{marking} for $X$ is an isomorphism of Chow motives
	  \[ \phi\colon\ 
	h(X)\stackrel{\cong}{\longrightarrow} M\ \ \hbox{in}\ \MM_{\rm rat}\ \] with $M$ a direct
	summand of a Chow motive of the form 
	$\oplus_{i} h(B_{i})(n_{i})$ cut out by an idempotent matrix $P \in
	\mathrm{End}(\oplus_i h(B_i)(n_i))$ whose entries are symmetrically
	distinguished cycles, where $B_{i}$ is an abelian variety and $n_{i}$ is an
	integer (the Tate twist). We refer to \cite[Definition~3.1]{FV} for the precise
	definition.
	
	Given a marking $\phi: h(X)\stackrel{\simeq}{\longrightarrow} M$, we define
	the subgroup of \emph{distinguished cycles} of $X$, denoted 
	$\operatorname{DA}^*_{\phi}(X)$, to be the pre-image of
	$\operatorname{DA}^*(M):= P_*\bigoplus_i \operatorname{DA}^{*-n_i}(B_i)$
	via the induced isomorphism
	$\phi_{*}:\A^*(X)\stackrel{\simeq}{\longrightarrow} \A^*(M)$. 
	
	Given another smooth projective variety $Y$ with a marking $\psi : h(Y) \to
	N$, the tensor product $\phi \otimes \psi : h(X\times Y) \to M\otimes N$
	naturally defines a marking for $X\times Y$. A morphism $f: X\to Y$ is said
	to be a \emph{distinguished morphism} if its graph is distinguished with respect
	to the product marking $\phi\otimes \psi$.

	The composition $$\operatorname{DA}^*_{\phi}(X)\hookrightarrow
	\A^*(X)\twoheadrightarrow \N^*(X)$$ is clearly bijective. In
	other words, $\phi$ provides a section (as graded vector spaces) of the natural
	projection $\A^*(X)\twoheadrightarrow\ \N^*(X)$. In \cite{FV}, 
	sufficient conditions on the marking $\phi$ are given such that
	$\operatorname{DA}^*_{\phi}(X)$ defines a $\QQ$-subalgebra of $\A^*(X)$\,:
	
	\begin{definition}[Definition 3.7 in \cite{FV}]\label{def:Star} 
		We say that the marking $\phi: h(X)\stackrel{\simeq}{\longrightarrow} M$
		satisfies the condition $(\star)$ if the following two conditions are
		satisfied\,: 
		\begin{itemize}
			\item[$(\star_{\mathrm{Mult}})$] the small diagonal $\delta_{X}$ belongs to
			$\operatorname{DA}^*_{\phi^{\otimes 3}}(X^{3})$\,; that is, under the
			induced
			isomorphism $\phi^{\otimes 3}_{*}:
			\A^*(X^{3})\stackrel{\simeq}{\longrightarrow} \A^*(M^{\otimes 3})$, the
			image
			of $\delta_{X}$ is symmetrically distinguished, \emph{i.e.} in
			$\operatorname{DA}^*(M^{\otimes 3})$.
			\item[$(\star_{\mathrm{Chern}})$]  all Chern classes $c_{i}(X)$ belong to
			$\operatorname{DA}^*_{\phi}(X)$.
		\end{itemize}
		If in addition $X$ is equipped with the action of a finite group $G$, we say
		that the  marking $\phi: h(X)\stackrel{\simeq}{\longrightarrow} M$  satisfies
		$(\star_G)$ if\,:
		\begin{itemize}
			\item[$(\star_G)$] the graph $g_X$ of $g: X\to X$ belongs to
			$\operatorname{DA}^*_{\phi^{\otimes 2}}(X^{2})$ for all $g\in G$. 
		\end{itemize}
	\end{definition}

The raison d'\^etre for condition $(\star)$ is its relation to the section property, as mentioned in the introduction:
		
	\begin{proposition}[Proposition 3.12 in \cite{FV}]\label{prop:Disting}
		If the marking $\phi: h(X)\stackrel{\simeq}{\longrightarrow} M$ satisfies the
		condition $(\star)$, then there is a section, as \emph{graded $\QQ$-algebras}, for
		the
		natural surjective morphism $\A^*(X)\to \N^*(X)$ such that all
		Chern
		classes of $X$ are in the image of this section.
		
		In other words, assuming $(\star)$ we have a graded $\QQ$-subalgebra
		$\operatorname{DA}_\phi^*(X)$ of the Chow ring $\A^*(X)$, which contains all
		the Chern
		classes of $X$ and is mapped isomorphically to $\N^*(X)$. Elements
		of $\operatorname{DA}_\phi^*(X)$ are called \emph{distinguished cycles}.
	\end{proposition}

	The raison d'\^etre for condition $(\star_G)$ is that it allows to easily treat quotients:
	
	\begin{proposition}\label{quot} Let $Y$ be a smooth projective variety verifying $(\star)$ and $(\star_G)$, for some finite group $G\subset\aut(Y)$.	
	Then $X/G$ has a marking such that $X/G$ verifies condition $(\star_{\mathrm{Mult}})$, and the quotient morphism $p\colon X\to X/G$ is distinguished. If 
	$p$ is \'etale then $X/G$ verifies condition $(\star)$.	
	\end{proposition}
	
\begin{proof} This is \cite[Proposition 4.12]{FV}.
\end{proof}

	We refer to \cite{FV} for 
	examples of varieties satisfying $(\star)$\,; for our purposes here, let us mention
	that these include abelian varieties, varieties with trivial Chow
	groups\footnote{A smooth projective variety $X$ is said to have \emph{trivial
			Chow groups} if $\A^\ast_{hom}(X)=0$, cf. \cite{V0}, \cite{Vo}.}, and hyperelliptic curves:
		
	\begin{proposition}\label{prop:markinghyperelliptic}
		Let $C$ be a hyperelliptic curve equipped with the action of the
		group
		$H\cong\ZZ/2\ZZ$ generated by the hyperelliptic involution. Then $C$ has a marking that satisfies the conditions $(\star)$ and $(\star_H)$,
		with the additional property that if $P$ is a fixed point of $H$, then the
		embedding $P\hookrightarrow C$ is distinguished.
	\end{proposition}
	
	\begin{proof} This is \cite[Proposition 3.3(\rom1)]{LV}.
	  \end{proof}

	The property $(\star)$ has great flexibility: as shown in \cite[Section 4]{FV}, 
 this property is stable under product, projectivization of vector bundles,
	and (under certain conditions) blow-ups.
	The relevant result for blow-ups is as follows:
	
	\begin{proposition}[\cite{FV}]\label{blowup}
		Let $X$ be a smooth projective variety and let $i:Y\hookrightarrow X$ be a
		closed smooth
		subvariety. Let $\tilde X$ be the blow-up of $X$ along $Y$ and let $E$ be the
		exceptional divisor, so that we have a cartesian diagram
		$$	\xymatrix{
			E\ar@{^{(}->}[r]^{j} \ar[d]_{p} & \tilde X\ar[d]^{\tau}\\
			Y \ar@{^{(}->}[r]^{i} & \ X\ .
		}$$
		If we have markings satisfying the condition $(\star)$ for $X$ and
		$Y$ such that $i:Y\hookrightarrow X$ is distinguished, then $E$ and $\tilde X$
		have  natural markings that satisfy $(\star)$ and are
		such that the morphisms $i,j,\tau$ and $p$ are all
		distinguished. 
		
		If in addition $X$ is equipped with the action of a finite group $G$ such that
		$G\cdot Y = Y$ and such that the markings of $X$ and $Y$ satisfy $(\star_G)$,
		then the natural markings of $E$ and $\tilde X$ also satisfy $(\star_G)$.\qed
	\end{proposition}
	
	\begin{proof} This is the content of \cite[Propositions 4.5 and 4.8]{FV}.
	\end{proof}

Let us also recall the following, which will come in useful in the proof of our main result (Theorem \ref{main}):
	
	\begin{proposition}\label{prop:crucial}
		Let $X$ be a smooth projective variety of dimension $n\geq 2$ with a 
		marking $\phi$
		satisfying the condition $(\star)$ of Definition \ref{def:Star}. Assume
		that
		the cohomology of $X$ is spanned by algebraic classes in degree $\neq n$. Then
		the graded $\QQ$-subalgebra $R^*(X) \subset \A^*(X)$ generated by divisors, Chern
		classes and by cycles that are the intersection of two cycles in $X$ of
		positive
		codimension (is contained in $\operatorname{DA}_\phi^*(X)$ and hence) injects into $\N^*(X)$.
	\end{proposition}      
	
\begin{proof} This is \cite[Proposition 2.10]{LV}.
\end{proof}

\section{Main result}

This section contains the proof of the main result, which is as follows:

\begin{theorem}\label{main} Let $X$ be a hyperelliptic Calabi--Yau $n$-fold. 

\noindent
(\rom1)
There exists a sequence of blow-ups $\wt{X}\to X$ such that $\wt{X}$ verifies condition $(\star)$, and hence for each $m\in\NN$ the $\QQ$-algebra epimorphism $\A^\ast(\wt{X}^m)\to \N^\ast(\wt{X}^m)$ admits a section, whose image contains the Chern classes of $\wt{X}^m$.

\noindent
(\rom2)
 Let $n\ge 2$. The $\QQ$-subalgebra $R^\ast(X)\subset \A^\ast(X)$ generated by divisors, Chern classes and by cycles that are intersections of two cycles of positive codimension injects into cohomology. In particular, the image of the intersection product
  \[ \A^i(X)\otimes \A^{j}(X)\ \ \to\ \A^{i+j}(X)\ \ \ \ (i,j>0) \]
  injects into cohomology.
  \end{theorem}

\begin{proof} 
To construct $\wt{X}$, we proceed as follows: starting from the singular double cover $\pi\colon\bar{X}\to\PP^n$, we first blow-up all $0$-dimensional loci $\pi^{-1}(H_{j_1}\cap\cdots\cap H_{j_n})$ (where the $j_i$ are pairwise distinct), then we blow-up all strict transforms of $1$-dimensional loci $\pi^{-1}(H_{j_1}\cap\cdots\cap H_{j_{n-1}})$ (where the $j_i$ are pairwise distinct), and so on (ending with codimension 2 loci).

This resolution process can be encoded in a cartesian diagram
  \begin{equation}\label{cart} \begin{array}[c]{cccccccc}
                \  \  \ \ \   \ Y_{n-1} &\xrightarrow{t_{n-2}}\ \ \ & \cdots & \cdots & \xrightarrow{}& Y_1 & \xrightarrow{t_0}& Y_0:= C^n \\
                      &&&&&&&\\
                  \ \   \ \  \ \  \ \ \ \ \downarrow{\scriptstyle p_{n-1}} \ \ \ &&&&& \ \ \ \  \downarrow{\scriptstyle p_1} &&   \downarrow{\scriptstyle p_0}\ \ \ \ \ \ \\
                      &&&&&&&\\
                      \wt{X}=:X_{n-1}\ \ \ \   & \xrightarrow{s_{n-2}}\ \ \ &      \cdots & \cdots & \xrightarrow{}& X_1 & \xrightarrow{s_0}& \ \ X_0:= \bar{X} \\       
                       &&&&&&&\\
                  \ \   \ \  \ \  \ \ \ \ \downarrow{\scriptstyle \pi_{n-1}} \ \ \ &&&&& \ \ \ \  \downarrow{\scriptstyle \pi_1} &&   \downarrow{\scriptstyle \pi_0}\ \ \ \ \ \ \\
                      &&&&&&&\\
                 \  \  \ \ \ \  \  P_{n-1}\ \ \ \   & \xrightarrow{r_{n-2}}\ \ \ &      \cdots & \cdots & \xrightarrow{}& P_1 & \xrightarrow{r_0}& \ \ P_0:= \PP^n \ .\\                          
                      \end{array}\end{equation}
                      Here, $r_0$ is the blow-up with center $\cup H_{j_1}\cap \cdots \cap H_{j_n}$, and
                       each $r_i$ is the blow-up with center $Q_i$ being the union of strict transforms of  $i$-dimensional intersections $H_{j_1}\cap \cdots\cap H_{j_{n-i}}$ (by construction, these strict transforms form a disjoint union). The arrows $s_i$ and $t_i$ are induced by $r_i$. Concretely, this means that $t_0$ is the blow-up with center
                          \[ Z_0:= {\displaystyle\bigcup_{\sigma\in\Sy_n}}\  \bigcup_{1\le j_1,\ldots,j_n\le 2n+2}\sigma\Bigl( q_{j_1}\times\cdots\times q_{j_n} \Bigr)\ \ \ \subset\  C^n\ ,\] 
  (where $q_1,\ldots,q_{2n+2}\in C$ are the Weierstrass points). Likewise, each $t_i$ is a blow-up with center
    \[ Z_i:= {\displaystyle\bigcup_{\sigma\in\Sy_n}}\  \bigcup_{1\le j_1,\ldots,j_{n-i}\le 2n+2}\overline{\sigma\bigl(  q_{j_1}\times \cdots\times q_{j_{n-i}} \bigr)} \ \ \ \subset\ Y_{i}\ \] 
     (where $\overline{a}$ means strict transform of $a$). By construction, $Z_i$ is a disjoint union of smooth irreducible components of dimension $i$. 
                     
                       The arrows $\pi_i$ are double covers. The arrow $p_0$ is the quotient morphism for the action of $G_0:=N\rtimes\Sy_n$, and the composition $\pi_0\circ p_0$ is the quotient morphism for the action of $H_0\cong (\ZZ/2\ZZ)^n \rtimes \Sy_n$.
                       Each arrow $p_i$ (resp. each composition $\pi_i\circ p_i$) is the quotient morphism for the action of the finite group $G_i$ (resp. $H_i$) on $Y_i$ obtained by lifting the action of $G_{i-1}$ (resp. $H_{i-1}$). 
%
     
 The idea is to prove property $(\star)$ for $\wt{X}$ inductively, moving from right to left in diagram \eqref{cart}. The induction base is $Y_0$:
 
 \begin{lemma}\label{y0} Let $G_0\subset H_0\subset \aut(C^n)$ be as above.
 The variety $Y_0:=C^n$ verifies conditions $(\star)$ and  $(\star_{G_0})$ and $(\star_{H_0})$.
  \end{lemma}
  
  \begin{proof}(of Lemma \ref{y0}.) The self-product $C^n$ verifies condition $(\star)$ because hyperelliptic curves verify $(\star)$ (Proposition \ref{prop:markinghyperelliptic}),
  and $(\star)$ is stable under products \cite[Proposition 4.1]{FV}. To check condition $(\star_{H_0})$ (which implies condition $(\star_{H_0})$), it suffices to check that the graph of $g$ is distinguished for any $g\in (\ZZ/2\ZZ)^n$ and for any $g\in\Sy_n$.
  For $g\in N$ this follows from the fact that the graph of the hyperelliptic involution is distinguished (Proposition \ref{prop:markinghyperelliptic}), plus the compatibility of group actions and products \cite[Proposition 4.1]{FV}.  For $g\in\Sy_n$, this follows from \cite[Remark 4.2]{FV}.  
   \end{proof}
   
  The induction step is as follows:
  
  \begin{lemma}\label{ind} Assume $Y_i$ verifies conditions $(\star)$
   and  $(\star_{G_i})$ and  $(\star_{H_i})$. 
  Then $Y_{i+1}$ verifies conditions $(\star)$.
   and  $(\star_{G_{i+1}})$ and $(\star_{H_{i+1}})$.  
  \end{lemma}
  
  \begin{proof}(of Lemma \ref{ind}.) This is an application of the general blow-up result Proposition \ref{blowup}. Let us check that all hypotheses of Proposition \ref{blowup} are met with. 
  The variety $Y_i$ and the center $Z_i$ of the blow-up verify $(\star)$: for $Y_i$ this is by assumption; for $Z_i$ this is true by induction, because $Z_i$ is $C^{i}$ blown-up along certain explicit centers (i.e. $Z_i$ is of the form $Y_i$ with a smaller value of $n$). 
 To see that the embedding $\iota_i\colon Z_i\hookrightarrow Y_i$ is distinguished, we note that its graph $\Gamma_{\iota_i}$ is the pullback of the graph of the embedding $\tau_i\colon Q_i\hookrightarrow P_i$ (here $Q_i$ and $P_i$ are as in the proof of Proposition \ref{p1}). The embedding $\tau_i$ is distinguished (indeed, one sees inductively that $P_i$ and $Q_i$ have trivial Chow groups, and so $\A^\ast(Q_i\times P_i)=\DA^\ast(Q_i\times P_i)$), and the assumption  $(\star_{H_i})$ implies that the quotient morphism
  $\pi_i\circ p_i$ is distinguished \cite[Proposition 4.12]{FV}. It follows that
  \[   \Gamma_{\iota_i}=  (\pi_i\circ p_i \times \pi_i\circ p_i)^\ast  \Gamma_{\tau_i}\ \ \ \in\ \DA^\ast(Z_i\times Y_i)\ .\]
  All hypotheses of Proposition \ref{blowup} being verified, this proves that $Y_{i+1}$ verifies condition $(\star)$.
  
  As for the group action, this follows from the second part of Proposition \ref{blowup}. Both $Y_i$ and the center $Z_i$ verify condition $(\star_{H_i})$ (and a fortiori $(\star_{G_i})$): for $Y_i$ this is by assumption, for $Z_i$ this is true by induction, because $Z_i$ is $C^{i}$ blown-up along certain explicit centers. 
  The second part of Proposition \ref{blowup} then guarantees that condition $(\star_{H_i})$ (and a fortiori $(\star_{G_i})$) carries over to $Y_{i+1}$.
    \end{proof}
  
The induction set up by Lemmas \ref{y0} and \ref{ind} yields that $Y_{n-1}$ verifies conditions $(\star)$ and  $(\star_{G_{n-1}})$. Using Proposition \ref{quot}, this implies that $\wt{X}= Y_{n-1}/G_{n-1}$ verifies condition $(\star_{\mathrm{Mult}})$. Let us now check condition  $(\star_{\mathrm{Chern}})$ for $\wt{X}$. For this, we view $\wt{X}$ as the double cover 
 \[\pi\colon \wt{X}\to \wt{P}:=P_{n-1}\] 
 branched along the smooth divisor $D\subset \wt{P}$ (obtained as strict transform of the hyperplane arrangement $\cup_{j=1}^{2n+2} H_j\subset\PP^n$ under the blow-ups $r_i$). The Chern classes of $\wt{X}$ can be expressed in terms of the Chern classes of $\wt{P}$ and the Chern classes of $\OO_{\wt{X}}(E)$ and $\OO_{\wt{X}}(2E)$, where we write $E\subset \wt{X}$ for the isomorphic pre-image of $D$ in $\wt{X}$
(cf. \cite[Proof of Proposition 4.12]{FV}). But there is equality
  \[  E= d\, (p_{n-1})_\ast (p_{n-1})^\ast\pi^\ast(D)\ \ \ \hbox{in}\ \A^1(\wt{X})\ \ \ \ (d\in\QQ)\ ,\] 
  and so (since $D\in \A^\ast(\wt{P})=\DA^\ast(\wt{P})$ and $\pi$ and $p_{n-1}$ are distinguished) one has $E\in \DA^\ast(\wt{X})$. Likewise, one has
  \[   \pi^\ast c_j(\wt{P})  =  d\, (p_{n-1})_\ast (p_{n-1})^\ast\pi^\ast c_j(\wt{P})\ \ \ \hbox{in}\ \A^\ast(\wt{X}) \ \ \ \ (d\in\QQ)\ ,\]
 and so also $\pi^\ast c_j(\wt{P})\in \DA^\ast(\wt{X})$.
We conclude that $c_j(\wt{X})\in  \DA^\ast(\wt{X})$, i.e. $\wt{X}$ verifies condition $(\star_{\mathrm{Chern}})$.

Since the Calabi--Yau variety $X$ was obtained from $\bar{X}$ by only blowing-up the codimension 2 loci $H_{i_1}\cap H_{i_2}$, there is a factorization
  \[  \wt{X}\ \xrightarrow{f}\ X\ \xrightarrow{}\ \bar{X}\ ,\]
  where both arrows are sequences of blow-ups with smooth centers. This proves (\rom1).
  
  We record the following:
  
  \begin{lemma}\label{compa} Let $\phi_i$ denote the marking for $X_i$ constructed above. Let $ S_i\subset X_i$ denote the center of the blow-up morphism $s_i\colon X_{i+1}\to X_i$, and let $\xi_i\colon S_i\to X_i$ denote the inclusion morphism. The graphs of $s_i$ and $\xi_i$ are distinguished (with respect to the markings $\phi_i$).
   \end{lemma}
   
   \begin{proof}(of Lemma \ref{compa}.)  There is a cartesian diagram
     \[   \begin{array}[c]{ccc}
                Y_{i+1} & \xrightarrow{t_i}& Y_i \\
                      &&\\
              \ \ \ \  \ \   \downarrow{\scriptstyle p_{i+1}} \ \ \ && \ \ \ \  \downarrow{\scriptstyle p_i} \\
                      &&\\
                      X_{i+1}  & \xrightarrow{s_{i}} &   \ \    X_i \ . \\       
                      \end{array} \]
                      
   The graph of $s_i\circ p_{i+1}$ equals the graph of $p_i\circ t_i$, which is distinguished by the construction above. Since the graph of $p_{i+1}$ is also distinguished (this is property $(\ast_{G_{i+1}})$), it follows that 
     \[ \Gamma_{s_i}=  {1\over \deg p_{i+1}}\ {}^t \Gamma_{p_{i+1}}  \circ \Gamma_{p_{i+1}}  \circ  \Gamma_{s_i} \]
     is distinguished.
     
     The argument for $\xi_i$ is similar, using the cartesian diagram
      \[   \begin{array}[c]{ccc}
                Z_{i} & \xrightarrow{\iota_i}& Y_i \\
                      &&\\
              \ \ \ \  \ \   \downarrow{} \ \ \ && \ \ \ \  \downarrow{\scriptstyle p_i} \\
                      &&\\
                      S_{i}  & \xrightarrow{\xi_{i}} &   \ \    X_i \ . \\       
                      \end{array} \]
                 \end{proof}
  
  \medskip
  \noindent
 (\rom2) Let us observe that $f\colon \wt{X}\to X$ is a sequence of blow-ups with smooth centers having trivial Chow groups, and so
  \begin{equation}\label{f} h(\wt{X})\cong  h(X)\oplus \bigoplus\one(\ast)\ \ \ \hbox{in}\ \MM_{\rm rat}\ .\end{equation}
  Using Proposition \ref{motive}(\rom2), we find that there is a decomposition
  \begin{equation}\label{g} h(\wt{X})\cong  M\oplus \bigoplus\one(\ast)\ \ \ \hbox{in}\ \MM_{\rm rat}\ ,\end{equation}
  where $M$ has cohomology concentrated in degree $n$. It then follows from Proposition \ref{prop:crucial} that the $\QQ$-subalgebra
   \[ R^\ast(\wt{X}):= \Bigl\langle  \A^1(\wt{X}), \, \A^i(\wt{X})\cdot \A^j(\wt{X}), \, c_k(\wt{X})\Bigr\rangle\ \ \ \subset\ \A^\ast(\wt{X})\ \ \ \ \ (i,j>0) \]
 is contained in $\DA^\ast(\wt{X})$ and hence injects into cohomology. Using pullback along the morphism $f\colon \wt{X}\to X$, this implies at once that the $\QQ$-subalgebra
   \[    \Bigl\langle  \A^1({X}), \, \A^i({X})\cdot \A^j({X})\Bigr\rangle\ \ \ \subset\ \A^\ast({X})\ \ \ \ \ (i,j>0) \] 
   also injects into cohomology. However, to get a statement that includes the Chern classes of $X$ some extra care is needed. Remark that one has
   \[  c_k(\wt{X}) =   f^\ast c_k(X) + R_k\ \ \ \hbox{in}\ \A^k(\wt{X})\ ,\]
   where $R_k$ is in the second part of the decomposition
   \begin{equation}\label{cho} \A^k(\wt{X})\cong \A^k(X)\oplus \QQ^s \end{equation}
   induced by \eqref{f}. Lemma \ref{compa} implies that the marking for $\wt{X}$ constructed above is induced (via the blow-up result Proposition \ref{blowup}) from a marking for $\bar{X}$, and so each trivial motive
   $\one(\ast)$ in \eqref{g} is marked by a trivial motive $\one(\ast)$. A fortiori, the same is true for \eqref{f}. This implies that the summand $\QQ^s$ in \eqref{cho} is in $\DA^\ast(\wt{X})$, and so $R_k\in \DA^\ast(\wt{X})$.\footnote{Alternatively, to see that  $R_k\in \DA^\ast(\wt{X})$ one could argue using Porteous' formula as in \cite[Proof of Proposition 4.4]{LV}.}
   But then one also has
   \[  f^\ast c_k(X)= c_k(\wt{X})-R_k\ \ \ \in\ \DA^\ast(\wt{X})\ ,\]
   and hence
   \[ f^\ast    \Bigl\langle  \A^1({X}), \, \A^i({X})\cdot \A^j({X}),\, c_k(X)\Bigr\rangle\ \ \ \subset\ \DA^\ast(\wt{X})\ \ \ \ \ (i,j>0) \ .\]    
  Since $\DA^\ast(\wt{X})$ injects into cohomology (under the cycle class map), we conclude that the $\QQ$-algebra
    \[ R^\ast(X):= \Bigl\langle  \A^1({X}), \, \A^i({X})\cdot \A^j({X}),\, c_k(X)\Bigr\rangle\ \ \ \ \ (i,j>0) \]
  also injects into cohomology (under the cycle class map).
    \end{proof}

\begin{remark} To prove our main result (Theorem \ref{main}(\rom2)), we show that for any hyperelliptic Calabi--Yau $n$-fold $X$ there exists a surjection $\wt{X}\to X$ such that $\wt{X}$ 
verifies condition $(\star)$.
This leaves open the question whether $X$ itself verifies condition $(\star)$ (which would provide an easier proof of the injectivity of Theorem \ref{main}(\rom2)).

Note that the hyperelliptic Calabi--Yau variety $X$ is a group quotient of some blow-up of $C^n$, with centers having trivial Chow groups. However, these centers do {\em not\/} behave well with respect to the group action (an irreducible codimension 2 center is not invariant under the action of $\Sy_n$), making it problematic to pass from the level of $Y_j$ to the level of the quotient $X_j$. For this reason we prove our main result in a slightly ``round-about'' way, blowing up some more in order to make the irreducible codimension 2 centers disjoint.
\end{remark}

\section{Further consequences}

\subsection{Voisin's conjecture}
	\label{svois}
	
	Voisin \cite{V9} has formulated the following intriguing conjecture, which is a
	particular instance of the Bloch--Beilinson conjectures.
	
	\begin{conjecture}[Voisin \cite{V9}]\label{conj} Let $X$ be a smooth projective
		variety of dimension $n$, with $h^{n,0}(X) = 1$ and $h^{j,0}(X)=0$
		for
		$0<j<n$. Then any two zero-cycles $a,a^\prime\in \A^n_{hom}(X)$ satisfy
		\[ a\times a^\prime = (-1)^n a^\prime\times a\ \ \ \hbox{in}\ \A^{2n}(X\times
		X)\ .\]
		(Here, $a\times a^\prime$ is the exterior product $(p_1)^\ast(a)\cdot
		(p_2)^\ast(a^\prime)\in \A^{2n}(X\times X)$, where $p_j$ is projection to the
		$j$-th factor.)
	\end{conjecture}

	For background and motivation for Conjecture \ref{conj}, \emph{cf.}
	\cite[Section 4.3.5.2]{Vo}.
	Conjecture \ref{conj} has been proven in some scattered special cases
	\cite{V9},
	\cite{moi}, \cite{desult}, \cite{tod}, \cite{24.5}, \cite{BLP}, \cite{LV}, \cite{Bur}, \cite{47}, but is still wide open for
	a general K3 surface.

	We now prove Voisin's conjecture for the Calabi--Yau varieties under consideration in the present paper\,: 
	
	\begin{proposition}\label{conjV} Let $X$ be a hyperelliptic Calabi--Yau variety of dimension $n$. Then Conjecture \ref{conj} is true for
		$X$:
		any $a,a^\prime\in \A^n_{hom}(X)$ satisfy
		\[ a\times a^\prime = (-1)^n \, a^\prime\times a\ \ \ \hbox{in}\
		\A^{2n}(X\times X)\ .\]
	\end{proposition}
	
	\begin{proof} 
According to Proposition \ref{motive}, we have a decomposition
		\[ h(X) = M\oplus \bigoplus {\mathds{1}}(\ast)\ \ \ \hbox{in}\ \MM_{\rm rat}, \]	
		with $H^j(M)=0$ for $j\not= n$, and $M$ isomorphic to a direct summand of $h(C^n/G)$, where $C$ is a hyperelliptic curve. Since the symmetric group $\Sy_n$ is contained in $G$, the motive $M$ is actually isomorphic to a direct summand of
		$h(C^{(n)})$. 
By Kimura finite-dimensionality \cite{Kim}, $M$ is then isomorphic to a direct summand of the
		motive
		\[  h^n(C^{(n)}):= (C^{(n)},\pi^n_{C^{(n)}},0)\ \ \ \in\ \MM_{\rm rat}\ ,\] 
		where $\pi^n_{C^{(n)}}$ is the Chow--K\"unneth
		projector defined by the choice of a point on $C$, i.e. $\pi^n_{C^{(n)}}$ corresponds to
		\[  \sum_{i_1+\cdots+i_n=n}\ \   \sum_{\sigma\in\Sy_n} \sigma\Bigl(   \pi_C^{i_1}\times\cdots \times \pi_C^{i_n} \Bigr) \ \ \in\ \A^n(C^n\times C^n)^{\Sy_n}\]
		under the natural isomorphism $\A^\ast(C^{(n)})\cong  \A^\ast(C^n)^{\Sy_n}$, and the $\pi_C^{i_j}$ are Chow--K\"unneth projectors of $C$.
				
	Let $\MM^\circ$ denote the category of birational motives \cite{KS}. There is a functor 
	\[ \bir\colon\  \MM_{\rm rat}^{\rm eff}\ \to\ \MM^\circ\ ,\] 
sending an effective motive $(X,p,0)$ to $(X,p\vert_{\A^n(X_{\C(X)})})$. This functor has the property that $\A_0(M)=\A_0(\bir(M))$ for any $M\in\MM_{\rm rat}^{\rm eff}$. Looking at the image of $	 h^n(C^{(n)})$ under $\bir$,
	one sees that all summands of $\pi^n_{C^{(n)}}$ where some $i_j$ is $0$ restrict to zero, 	i.e. one has
	\[  \bir \bigl( h^n(C^{(n)})\bigr)=\bir \bigl(\sym^n h^1(C)\bigr)\ \ \ \hbox{in}\ \MM^\circ\ .\]
	Writing $J:=\jac(C)$ for the Jacobian of $C$ and recalling that there is an isomorphism $h^1(C)\cong h^1(J)$, it follows that there are isomorphisms of birational motives
	\[  \bir \bigl( h^n(C^{(n)})\bigr)\cong\bir \bigl(\sym^n h^1(C)\bigr)\cong  \bir \bigl(\sym^n h^1(J)\bigr)\cong \bir\bigl( h^n(J)\bigr)   \ \ \ \hbox{in}\ \MM^\circ\ ,\]
	where $h^\ast(J)$ refers to the Deninger--Murre Chow--K\"unneth decomposition for abelian schemes \cite{DM} (for the properties of Chow motives of abelian varieties that we use here, cf. \cite[Section 5]{Sch}). In particular, taking Chow groups we find a split injection
	\[ \begin{split}  \Gamma_\ast\colon\ \  \A_0(M)\ \hookrightarrow\ \A_0 (  h^n(C^{(n)})&=\A_0\bigl(   \bir ( h^n(C^{(n)})) \bigr)\\ &=  \A_0\bigl(\bir( h^n(J))\bigr)= \A_0(h^n(J))=\A^n_{(n)}(J)\ ,\\
	\end{split} \]
	where $\A^\ast_{(\ast)}(J)$ refers to Beauville's eigenspace decomposition \cite{Beau}. By the same argument, there is also a split injection
	\[  (\Gamma\times\Gamma)_\ast\colon\ \   \A_0(M\otimes M)\ \hookrightarrow\ \A_0(J\times J)  \ . \]
	This fits into a commutative diagram
	\[ \begin{array}[c]{ccc} \A^n(M)\otimes \A^n(M) & \xrightarrow{(\Gamma_\ast,\Gamma_\ast)}& \A^n_{(n)}(J)\otimes \A^n_{(n)}(J)\\
	                                               &&\\
	                                               \ \ \downarrow{\scriptstyle \Phi}&&   \ \ \downarrow{\scriptstyle \Phi}\\
	                                               &&\\
	                                           \A^{2n}(M\otimes M) & \xrightarrow{(\Gamma\times\Gamma)_\ast}&\ \A^{2n}_{}(J\times J)\ ,\\
	                                              \end{array} \]
	                                              where $\Phi$ sends $(a,a^\prime)$ to $a\times a^\prime - (-1)^n \, a^\prime\times a$. We are thus reduced to a general statement about zero-cycles on abelian varieties:

		\begin{proposition}[Voisin \cite{Vo}]\label{abvar} 
			Let $B$ be an abelian variety of dimension $n$, and $a,a^\prime\in
			\A^n_{(n)}(B)$. Then
			\[a\times a^\prime=(-1)^n \, a^\prime\times a\ \ \ \hbox{in}\
			\A^{2n}(B\times B) .\]
		\end{proposition}
		
		\begin{proof} This is \cite[Example 4.40]{Vo}. A generalization (and an alternative proof) is given in \cite[Theorem 4.1]{Ch}.
		\end{proof} 
		
		This concludes the proof of the theorem.
	\end{proof}

	\subsection{Voevodsky's conjecture}
	\label{svoe}
	
	This subsection contains an application of our results to Voevodsky's
	conjecture on smash-equivalence.

	\begin{definition}[Voevodsky \cite{Voe}]\label{sm} Let $X$ be a smooth
		projective variety. A cycle $a\in \A^i(X)$ is called {\em smash-nilpotent\/} 
		if there exists $m\in\NN$ such that
		\[ \begin{array}[c]{ccc}  a^m:= &\underbrace{a\times\cdots\times a}_{(m\hbox{
				times})}&=0\ \ \hbox{in}\ 
		\A^{mi}(X\times\cdots\times
		X)_{}\ .
		\end{array}\]
		\vskip0.6cm
		
		Two cycles $a,a^\prime$ are called {\em smash-equivalent\/} if their
		difference
		$a-a^\prime$ is smash-nilpotent. We will write $\A^i_\otimes(X)\subseteq
		\A^i(X)$ for the subgroup of smash-nilpotent cycles.
	\end{definition}
	
	\begin{conjecture}[Voevodsky \cite{Voe}]\label{voe} Let $X$ be a smooth
		projective variety. Then
		\[  \A^i_{hom}(X)\ \subseteq\ \A^i_\otimes(X)\ \ \ \hbox{for\ all\ }i\ .\]
	\end{conjecture}
	
	\begin{remark} It is known \cite[Th\'eor\`eme 3.33]{An} that Conjecture
		\ref{voe} for all smooth projective varieties implies (and is strictly
		stronger
		than) Kimura's conjecture ``all smooth projective varieties have
		finite-dimensional motive'' \cite{Kim}.
	\end{remark}

	Let us now verify Voevodsky's conjecture for hyperelliptic Calabi--Yau varieties of odd dimension:
	
	\begin{theorem}\label{conjVoe} Let $X$ be a hyperelliptic Calabi--Yau variety. Assume that $n:=\dim X$ is odd.
		Then
		\[  \A^i_{hom}(X)\ \subseteq\ \A^i_\otimes(X)\ \ \ \hbox{for\ all\ }i.\]
	\end{theorem}
	
	\begin{proof} 
		According to Proposition \ref{motive}, we have a decomposition
		\[ h(X) = M\oplus \bigoplus {\mathds{1}}(\ast)\ \ \ \hbox{in}\ \MM_{\rm rat}, \]	
		with $H^j(M)=0$ for $j\not= n$, and $M$ isomorphic to a direct summand of
		$h(C^n)$. 				
		By Kimura finite-dimensionality, $M$ is isomorphic to a direct summand of the
		motive
		$(C^n,\pi^n,0)$, where $\pi^n$ is any Chow--K\"unneth
		projector on the degree-$n$ cohomology.		
		But the Chow motive $(C^n,\pi^n,0)$ is oddly
		finite-dimensional (in the sense of \cite{Kim}).
		Hence, together with the fact that $\A^{i}_{hom}({X})=
		\A^i_{hom}(M)$, the theorem is implied by the fact that 
		\[ \A^\ast_{}(M)
		\ \subseteq\ \A^\ast_\otimes(M)\] for any oddly finite-dimensional
		Chow motive $M$  (this is due to Kimura
		\cite[Proposition 6.1]{Kim}, and is also used in \cite{KSeb}). 
	\end{proof}

\section{Some open questions}
\label{quest}

\begin{question} Easy examples of Calabi--Yau varieties $Y$ are given by smooth complete intersections of $n+1$ quadrics in $\PP^{2n+1}$. An interesting special case is when $Y$ is defined by equations of the form
    \[   \begin{cases}  x_0^2+x_1^2+\cdots \  \cdots \ \cdots+ x_{2n+1}^2&=0\ ,\\
                          \lambda_0 x_0^2+\lambda_1 x_1^2+\cdots + \lambda_{2n+1}x_{2n+1}^2&=0\ ,\\
                            \lambda_0^2 x_0^2+\lambda_1^2 x_1^2+\cdots + \lambda_{2n+1}^2 x_{2n+1}^2&=0\ ,\\
                          \ \ \ \vdots\\
                         \ \ \  \vdots\\
                           \lambda_0^n x_0^2+\lambda_1^n x_1^2+\cdots + \lambda_{2n+1}^n x_{2n+1}^2&=0\ ,\\     
                           \end{cases}\]
                           where $\lambda_0,\ldots,\lambda_{2n+1}\in\C$ are distinct numbers.
                       Such varieties $Y$ are isomorphic to a quotient $D^n/G$, where $D$ is a curve and $G\subset\aut(D^n)$ a finite group \cite[Theorem 2.4.2]{Tera}.
                       Moreover, $Y$ is related to the double cover $\bar{X}$ of Theorem \ref{main}: one has $\bar{X}\cong Y/H$ for some finite group $H$ \cite[Proposition 2.4.4]{Tera}, \cite[Section 2.2]{GSSZ}. 
                       
                       Can one prove condition $(\star)$ for $Y$ ? The problem is that the curve $D$ is {\em not\/} hyperelliptic; $D$ is a finite \'etale cover of a hyperelliptic curve, and it is not clear whether $D$ verifies condition $(\star)$ (this is closely related to the fact that as far as I am aware the only curves known to have a multiplicative Chow--K\"unneth decomposition, in the sense of \cite{SV}, are hyperelliptic curves).
\end{question}

\begin{question} Let $Y$ be a Calabi--Yau variety as in Proposition \ref{p1}, i.e. arising from a general arrangement of hyperplanes (without the condition that the hyperplanes osculate a rational normal curve). Is it still true that the image of intersection product
   \[ \A^i(Y)\otimes \A^{j}(Y)\ \ \to\ \A^{i+j}(Y)\ \ \ \ (i,j>0) \]
  injects into cohomology ? 
 The problem is that for these $Y$, it is not even known that they have motive of abelian type, so that one cannot benefit from the formalism of \cite{FV}.
\end{question}

 \vskip1cm
\begin{nonumberingt} Thanks to the referee for very helpful comments. Thanks to Kai and Len for enjoying Buurman en Buurman as much as I do.
\end{nonumberingt}


\vskip1cm
\begin{thebibliography}{dlPG99}

\bibitem{An0} Y. Andr\'e, Une introduction aux motifs (motifs purs, motifs mixtes, p\'eriodes), Panoramas et Synth\`eses [Panoramas and Syntheses], vol. 17, Soci\'et\'e Math\'ematique de France, Paris 2004,

\bibitem{An} Y. Andr\'e, Motifs de dimension finie (d'apr\`es S.-I. Kimura, P. O'Sullivan,...), S\'eminaire Bourbaki 2003/2004, Ast\'erisque 299 Exp. No. 929, viii, 115---145,
			
\bibitem{Beau} A. Beauville, Sur l'anneau de Chow d'une variété abélienne. Math. Ann. 273 (1986), 647---651,
		
\bibitem{Beau3} A. Beauville, On the splitting of the Bloch--Beilinson filtration, in: Algebraic cycles and motives (J. Nagel and C. Peters, editors), London Math. Soc. Lecture Notes 344, Cambridge University Press 2007,

\bibitem{BV} A. Beauville and C. Voisin, On the Chow ring of a K3 surface, J. Alg. Geom. 13 (2004), 417---426,

\bibitem{BLP} G. Bini, R. Laterveer and G. Pacienza, Voisin's conjecture for zero-cycles on Calabi--Yau varieties and their mirrors, Advances in Geometry 20 no. 1 (2020), 91---108,



\bibitem{Bur} D. Burek, Higher-dimensional Calabi--Yau manifolds of Kummer type, Math. Nachrichten 293 no. 4 (2020), 638---650,

\bibitem{CH} S. Cynk and K. Hulek, Construction and examples of higher-dimensional modular Calabi--Yau manifolds, Can. Math. Bull. 50 no. 4 (2007), 486---503,

\bibitem{CK} S. Cynk and B. Kocel--Cynk, Classification of double octic Calabi--Yau threefolds with $h^{1,2}\le 1$ defined by an arrangement of eight planes, Commun. Contemp. Math. 
22 no. 1 (2020),

\bibitem{CM} S. Cynk and C. Meyer, Geometry and arithmetic of certain double octic Calabi--Yau threefolds, Can. Math. Bull. 48 no. 2 (2005), 180---194,

\bibitem{CM2} S. Cynk and C. Meyer, Modularity of some non-rigid double octic Calabi--Yau threefolds, Rocky Mountain J. Math. 38 no. 6 (2008), 1937---1958,

\bibitem{CSS} S. Cynk, M. Sch\"utt and D. van Straten, Hilbert modularity of some double octic Calabi--Yau threefolds. J. Number Theory 210 (2020), 313---332,

\bibitem{CS} S. Cynk and D. van Straten, Periods of rigid double octic Calabi--Yau threefolds, Ann. Pol. Math. 123 Part 1 (2019), 243---258,

\bibitem{DM}  Ch. Deninger and J. Murre, Motivic decomposition of abelian schemes and the Fourier transform, J. Reine Angew. Math. 422 (1991), 201---219,
		
\bibitem{Dol} I. Dolgachev, Weighted projective varieties, in: Group actions and vector fields (Vancouver, B.C., 1981) ed. J. B. Carrell, Lecture Notes in Math., vol. 956, Springer, Berlin, 1982, 34---71,
%
%
%
%

\bibitem{Fu} L. Fu, Decomposition of small diagonals and Chow rings of hypersurfaces and Calabi--Yau complete intersections, Advances in Mathematics 244 (2013), 894---924,

\bibitem{FV} L. Fu and Ch. Vial, Distinguished cycles on varieties with motive of abelian type and the section property, J. Alg. Geom. 29 (2020), 53---107,

\bibitem{F} W. Fulton, Intersection theory, Springer--Verlag Ergebnisse der Mathematik, Berlin Heidelberg New York Tokyo 1984,

\bibitem{GSSZ} R. Gerkmann, M. Sheng, D. van Straten and K. Zuo, On the monodromy of the moduli space of Calabi--Yau threefolds coming from eight planes in $\PP^3$, Math. Ann. 
355 (2013), 187---214,

\bibitem{IL} C. Ingalls and A. Logan, On the Cynk--Hulek criterion for crepant resolutions of double covers, arXiv:2006.14981v2,

%

\bibitem{KSeb} B. Kahn and R. Sebastian, Smash-nilpotent cycles on abelian 3-folds, Math. Res. Letters 16 (2009), 1007---1010,

\bibitem{KS} B. Kahn and R. Sujatha, Birational motives I: Pure birational motives, Ann. K-Theory 1 no. 4 (2016), 379---440,

\bibitem{Kim} S.-I.~Kimura, Chow groups are finite dimensional, in some sense, Math. Ann. 331 no 1 (2005), 173---201, 

\bibitem{moi} R. Laterveer, Some results on a conjecture of Voisin for surfaces of geometric genus one, Boll. Unione Mat. Italiana 9 no. 4 (2016), 435---452,
		
\bibitem{desult} R. Laterveer, Some desultory remarks concerning algebraic cycles and Calabi--Yau threefolds, Rend. Circ. Mat. Palermo 65 (2) (2016), 333---344,
		
\bibitem{tod} R. Laterveer, Algebraic cycles and Todorov surfaces, Kyoto Journal of Mathematics 58 no. 3 (2018), 493---527,
		
\bibitem{24.5} R. Laterveer, Some Calabi--Yau fourfolds verifying Voisin's conjecture, Ricerche di Matematica 67 no. 2 (2018), 401---411,
		
\bibitem{47} R. Laterveer,  Zero-cycles on self-products of varieties: some elementary examples verifying Voisin's conjecture, Bolletino Unione Mat. Italiana,
doi: 10.1007/s40574-020-00259-0,		
		
\bibitem{LV} R. Laterveer and Ch. Vial, On the Chow ring of Cynk--Hulek Calabi--Yau varieties and Schreieder varieties, Canadian Journal of Math. 72 no 2 (2020), 505---536,

\bibitem{MS} K. Matsumoto, T. Sasaki and M. Yoshida, The monodromy of the period map of a 4-parameter family
of K3 surfaces and the hypergeometric function of type $(3,6)$, Int. J. Math. 3 no. 1 (1992), 1---164,

\bibitem{MT} K. Matsumoto and T. Terasoma, Arithmetic-geometric means for hyperelliptic curves and Calabi-Yau varieties, Int. J. Math. 21 no. 7 (2010), 939---949,

\bibitem{Mey} C. Meyer, Modular Calabi--Yau threefolds. Fields Institute Monographs, 22. American Mathematical Society, Providence RI 2005,

\bibitem{MNP} J. Murre, J. Nagel and C. Peters, Lectures on the theory of pure motives, Amer. Math. Soc. University Lecture Series 61, Providence 2013,


\bibitem{OS} P. O'Sullivan, Algebraic cycles on an abelian variety, J. Reine Angew. Math. 654 (2011), 1---81,
		
\bibitem{Par} K. Paranjape, Abelian varieties associated to certain K3 surfaces, Comp. Math. 68 no. 1 (1988), 11---22,		

\bibitem{Sch} T. Scholl, Classical motives, in: Motives (U. Jannsen et alii, eds.), Proceedings of Symposia in Pure Mathematics Vol. 55 (1994), Part 1, 



\bibitem{Tera} T. Terasoma, Complete intersections of hypersurfaces---the Fermat case and the quadric case, Japan J. Math. (N.S.) 14 no. 2 (1988), 309---384,

\bibitem{Tera2} T. Terasoma, Algebraic correspondences between genus three curves and certain Calabi--Yau varieties, American Journal of Mathematics 132 no. 1 (2010), 181---200,

\bibitem{Ch} Ch. Vial, Generic cycles, Lefschetz representations and the generalized Hodge and Bloch conjectures for abelian varieties, Annali della Scuola Normale Superiore di Pisa 21 (2020), 1389---1429,

\bibitem{Voe} V. Voevodsky, A nilpotence theorem for cycles algebraically equivalent to zero, Internat. Math. Research Notices 4 (1995), 187---198,
		
\bibitem{V9} C. Voisin, Remarks on zero-cycles of self-products of varieties, in: Moduli of vector bundles, Proceedings of the Taniguchi Congress  (M. Maruyama, ed.), Marcel Dekker New York Basel Hong Kong 1994,
		
\bibitem{V0} C. Voisin, The generalized Hodge and Bloch conjectures are equivalent for general complete intersections, Ann. Sci. Ecole Norm. Sup. 46, fascicule 3 (2013), 449---475,
		
\bibitem{Vo} C. Voisin, Chow Rings, Decomposition of the Diagonal, and the Topology of Families, Princeton University Press, Princeton and Oxford, 2014,
	
%
\bibitem{Yos} M. Yoshida, Hypergeometric Functions, My Love: Modular Interpretations of Configuration Spaces,
Aspects of Mathematics, E32. Friedr. Vieweg \& Sohn, Braunschweig 1997,



\end{thebibliography}
\end{document}